\documentclass[11pt]{article}
\renewcommand{\baselinestretch}{1.2}
\newcommand{\dated}{\mbox{} \hfill {\small [{\tt \today}]}} % theorems etc.
\newtheorem{theorem}{Theorem}[section]
\newtheorem{lemma}[theorem]{Lemma}
\newtheorem{corollary}[theorem]{Corollary}
\newtheorem{proposition}[theorem]{Proposition}
\newtheorem{df}[theorem]{Definition}
\newenvironment{definition}{\begin{df} \rm}{\end{df}}

 \usepackage{amsmath,amssymb,amscd}
\newcommand{\pf}[1]{\trivlist \item[\hskip\labelsep\it #1\ ]}
\newcommand{\varpf}[1]{\trivlist \item[\hskip\labelsep\sc #1:]}
\newcommand{\qedbox}{$\rlap{$\sqcap$}\sqcup$}
\newcommand{\qed}{\qquad \qedbox \endtrivlist}
\newcommand{\varqed}{\hfill \rule{0.6em}{0.6em} \endtrivlist}
\newenvironment{proof}{\pf{Proof}}{\qed}

\newenvironment{remark}{\pf{Remark}}{\endtrivlist}
\newenvironment{remarks}{\pf{Remarks} 
   \begin{enumerate}}{\end{enumerate} \endtrivlist}
\newenvironment{example}{\pf{Example}}{\endtrivlist}
\newenvironment{examples}{\pf{Examples} 
   \begin{enumerate}}{\end{enumerate} \endtrivlist}
\newenvironment{items}{
  \begin{enumerate} 
                    
  }{\end{enumerate}}

\newenvironment{keywords}{\noindent\small {\it Keywords\/}:}{\vskip 4pt}
\newenvironment{classification}{\noindent\small 2000 {\it Mathematics Subject
Classification\/}:}{\vskip 12pt}

\newcommand{\comps}{{\mathbb C}}

\newcommand{\posints}{{\mathbb N}}

\newcommand{\tensor}{\otimes}

\newcommand{\Tensor}{\hat{\otimes}}

\newcommand{\cstar}{{C^\ast}}

\newcommand{\id}{{\mathrm{id}}}
\newcommand{\cb}{{\mathrm{cb}}}

\newcommand{\A}{{\mathfrak A}}

\newcommand{\M}{{\mathfrak M}}

\newcommand{\CB}{{\cal CB}}

\newcommand{\VN}{\operatorname{VN}}

\newcommand{\AM}{\mathrm{AM}}
\title{The amenability constant of the Fourier algebra}
\author{{\it Volker Runde}\thanks{Research supported by NSERC under grant no.\ 227043-04.}}
\date{}
\begin{document}
\maketitle
\begin{abstract}
For a locally compact group $G$, let $A(G)$ denote its Fourier algebra and
$\hat{G}$ its dual object, i.e.\ the collection of equivalence classes of
unitary representations of $G$. We show that the amenability constant of $A(G)$
is less than or equal to $\sup \{ \deg(\pi) : \pi \in \hat{G} \}$ and that it
is equal to one if and only if $G$ is abelian.
\end{abstract}
\section*{Introduction}
The theory of amenable Banach algebras begins with B.\ E.\ Johnson's memoir
\cite{Joh1}. The choice of terminology is motivated by 
\cite[Theorem 2.5]{Joh1}: a locally compact group is amenable (in the usual
sense; see \cite{Pie}, for example), if and only if
its group algebra $L^1(G)$ is an amenable Banach algebra. For a modern
account of the theory of amenable Banach algebras, see \cite{LoA}.
\par
The Fourier algebra $A(G)$ of an arbitrary locally compact group $G$ was
introduced by P.\ Eymard in \cite{Eym}. If $G$ is abelian, then the Fourier
transform yields an isometric isomorphism of $A(G)$ and $L^1(\hat{G})$,
where $\hat{G}$ is the dual group of $G$. (In the framework of Kac algebras,
this extends to a duality between $L^1(G)$ and $A(G)$ for arbitrary $G$;
see \cite{ES}.) Since amenable Banach algebras have bounded approximate 
identities, Leptin's theorem (\cite{Lep}) yields immediately that $A(G)$ can 
be amenable only if $G$ is amenable.
\par
Nevertheless, the tempting conjecture that a locally compact group $G$ is 
amenable if and {\it only if\/} $A(G)$ is amenable, turned out to be wrong, as
Johnson showed in \cite{Joh3}. For any locally compact group $G$, let 
$\hat{G}$ denote its {\it dual object\/}, i.e.\ the collection of all 
equivalence classes of (continuous) irreducible unitary representations of 
$G$. For $\pi \in \hat{G}$, let let $\deg(\pi)$ denote its {\it degree\/}, 
i.e.\ the dimension of the corresponding Hilbert space. For compact $G$, 
Johnson showed: If $G$ is infinite such that such that 
$\{ \pi \in \hat{G} : \deg(\pi) = n \}$ is finite for each $n \in \posints$, 
the Fourier algebra cannot be amenable. Hence, for example, 
$A(\operatorname{SO}(3))$ is not amenable. 
\par
This leaves the problem to characterize those locally compact groups $G$ for
which $A(G)$ is amenable (\cite[Problem 14]{LoA}). On the positive side,
$A(G) \cong L^1(\hat{G})$ is amenable whenever $G$ is abelian, and $A(G)$
is trivially amenable if $G$ is finite. With a little more effort, one can
show that, if $G$ is {\it almost abelian\/}, i.e.\ has an abelian subgroup
of finite index, then $A(G)$ is still amenable (\cite[Theorem 4.1]{LLW}).
Eventually, the locally compact groups $G$ with an amenable Fourier algebra 
were characterized by B.\ E.\ Forrest and the author: $A(G)$ is amenable if
and {\it only if\/} $G$ is almost abelian (\cite[Theorem 2.3]{FR}).
\par
In the present note, we will pick up another line of investigation begun in
\cite{Joh3} (and continued in \cite{LLW}). Suppose that $A(G)$ is amenable, 
so that it makes sense to speak 
of its amenability constant, which we denote by $\AM_{A(G)}$. For finite
$G$, Johnson, in \cite{Joh3}, derived a remarkable formula that allows to 
compute $\AM_{A(G)}$ in terms of the degrees of the irreducible
unitary representations of $G$, namely
\begin{equation} \label{Barry}
  \AM_{A(G)} = 
  \frac{\sum_{\pi \in \hat{G}} \deg(\pi)^3}{\sum_{\pi \in \hat{G}} \deg(\pi)^2}.
\end{equation}
\par
From (\ref{Barry}), it is immediate that the following are true for finite
$G$:
\begin{itemize}
\item $\AM_{A(G)} \leq \deg(G) := \sup \{ \deg(\pi) : \pi \in \hat{G} \}$;
\item $\AM_{A(G)} = 1$ if and only if $G$ is abelian.
\end{itemize}
\par
It is the purpose of this note to show that these two statements on 
$\AM_{A(G)}$ are true not only if $G$ is finite, but for all locally compact
groups $G$. As a by-product, we obtain an alternative approach to 
\cite[Theorem 2.3]{FR}.
\section{Amenability preliminaries}
Johnson's original definition of an amenable Banach algebra was in terms of
cohomology groups (\cite{Joh1}). We prefer to give another approach, which is 
based on a characterization of amenable Banach algebras from \cite{Joh2}.
\par
Following \cite{ER}, we denote the (completed) Banach space tensor product
by $\tensor^\gamma$. If $\A$ is a Banach algebra, then $\A \tensor^\gamma \A$
becomes a  Banach $\A$-bimodule via
\[
  a \cdot (x \tensor y) := ax \tensor y \quad\text{and}\quad (x \tensor y) \cdot a := x \tensor ya \qquad (a,x,y \in \A).
\]
The product of $\A$ induces a homomorphism $\Delta_\A \!: \A \tensor^\gamma \A 
\to \A$ of Banach $\A$-bimodules.
\begin{definition}
Let $\A$ be a Banach algebra. An {\it approximate diagonal\/} for $\A$ is
a bounded net $( d_\alpha )_\alpha$ in $\A \tensor^\gamma \A$ such that
\begin{equation} \label{diag1}
  a \cdot d_\alpha - d_\alpha \cdot a \to 0 \qquad (a \in \A)
\end{equation}
and
\begin{equation} \label{diag2}
  a \Delta_\A d_\alpha \to a \qquad (a \in \A).
\end{equation}
\end{definition}
\par
By \cite{Joh2}, a Banach algebra is amenable if and only if it has an
approximate diagonal. The advantage of using approximate diagonals to 
define  amenable Banach algebras is that approximate diagonals allow to
introduce a quantitative aspect into the notion of amenability:
\begin{definition}
A Banach algebra $\A$ is called $C$-amenable with $C \geq 0$ if there is
an approximate diagonal for $\A$ bounded by $C$.
\end{definition}
\begin{remarks}
\item In view of \cite{Joh2}, a Banach algebra is amenable if and only if
it is $C$-amenable for some $C \geq 0$.
\item By (\ref{diag2}) it is impossible for any Banach algebra to be
$C$-amenable with $C < 1$.
\end{remarks}
\begin{definition} \label{AMdef}
Let $\A$ be a Banach algebra. The {\it amenability constant\/} of $\A$ is 
defined as
\[
  \AM_\A := \inf \{ C \geq 0 : \text{$\A$ is $C$-amenable} \}.
\]
\end{definition}
\begin{remarks}
\item In terms of Definition \ref{AMdef}, $\A$ is amenable if and only if
$\AM_\A < \infty$.
\item The infimum in Definition \ref{AMdef} is easily seen to be a minimum.
\end{remarks}
\begin{examples}
\item Let $G$ be a locally compact group. Then $G$ is amenable, if and only if
$L^1(G)$ is $1$-amenable (\cite[Corollary 1.11]{Sto}).
Hence, we either have $\AM_{L^1(G)} = \infty$ or $\AM_{L^1(G)} =1$ depending on
whether $G$ is amenable or not.
\item Let $\A$ be a $\cstar$-algebra. Then $\A$ is amenable if and only if it
is nuclear (see \cite[Chapter 6]{LoA} for a self-contained exposition).
By \cite[Theorem 3.1]{Haa}, if $\A$ is nuclear, then it is already $1$-amenable. 
We thus have again a dichotomy that either $\AM_\A = \infty$ of $\AM_\A = 1$.
\item Let $G$ be a finite group. Then $\AM_{A(G)}$ can be explicitly computed
through (\ref{Barry}). From (\ref{Barry}), it follows immediately that 
$\AM_{A(G)} = 1$ if and only if $G$ is abelian, but more is true: if $G$ is 
not abelian, then $\AM_{A(G)} \geq \frac{3}{2}$ must hold
(\cite[Proposition 4.3]{Joh3}).
Another consequence of (\ref{Barry}) is that, if $H$ is another finite group,
we have
\[
  \AM_{A(G \times H)} = \AM_{A(G)} \, \AM_{A(H)}.
\]
Consequently, if $G$ is not abelian, $\AM_{A(G^n)} \geq 
\left( \frac{3}{2} \right)^n$ can be arbitrarily large.
\end{examples}
\section{An estimate from above for $\AM_{A(G)}$}
In this section, we shall extend (\ref{Barry}) to general locally compact
groups in the sense that we shall show, that for any locally compact group
$G$, the inequality $\AM_{A(G)} \leq \deg(G)$ holds.
\par
We require some background from the theory of operator spaces, for which we
refer to \cite{ER}, whose notation we adopt. In particular, for a linear space
$E$ and $n \in \posints$, the symbol $M_n(E)$ stands for the 
$n \times n$-matrices with entries from $E$, and if $F$ is another linear space, 
and $T \!: E \to F$ is linear, then the $n$-th amplification of $T$ --- 
from $M_n(E)$ to $M_n(F)$ --- is denoted by $T_n$.
\par
Our first lemma, is a minor generalization of \cite[Proposition 2.2.6]{ER} and 
has an almost identical proof:
\begin{lemma} \label{smithie}
Let $E$ be an operator space, let $\A$ be a commutative $\cstar$-algebra,
and let $n \in \posints$. Then every bounded linear map $T \!: E \to
M_n(\A)$ is completely bounded such that $\| T \|_\cb = \| T_n \|$.
\end{lemma}
\begin{proof}
Let $\Omega$ be a locally compact Hausdorff space such that $\A \cong
{\cal C}_0(\Omega)$. We may identify $M_n(\A)$ with ${\cal C}_0(\Omega,M_n)$.
For $\omega \in \Omega$, let 
\[
  T^\omega \!: E \to M_n, \quad x \mapsto (Tx)(\omega).
\]
By Smith's lemma (\cite[Proposition 2.2.2]{ER}), each map 
$T^\omega$ is completely bounded with $\| T^\omega \|_\cb = \| T^\omega_n \|$, so that
\begin{equation} \label{roger}
   \| T^\omega \|_\cb = \| T^\omega_m \| = \| T^\omega_n \|
   \qquad
   (m \in \posints, \, m \geq n).
\end{equation}
\par
Let $m \in \posints$ with $m \geq n$. Then we have:
\[
  \begin{split}
  \| T_m \| & = \sup 
  \left\{ \| T_m x \|_{{\cal C}_0(\Omega,M_{mn})} : 
  x \in M_m(E), \, \| x \|_{M_m(E)} \leq 1 \right\}
  \\
  & =  \sup 
  \left\{ \| T_m^\omega x \|_{M_{mn}} : \omega \in \Omega, \, 
  x \in M_m(E), \, \| x \|_{M_m(E)} \leq 1 \right\} \\
  & = \sup \{ \| T_m^\omega \| : \omega \in \Omega \} \\
  & = \sup \{ \| T_n^\omega \| : \omega \in \Omega \}, \qquad\text{by
  (\ref{roger})}, \\
  & =  \sup 
  \left\{ \| T_n^\omega x \|_{M_{n^2}} : \omega \in \Omega, \, 
  x \in M_n(E), \, \| x \|_{M_n(E)} \leq 1 \right\} \\
  & = \sup 
  \left\{ \| T_n x \|_{{\cal C}_0\left(\Omega,M_{n^2}\right)} : 
  x \in M_n(E), \, \| x \|_{M_n(E)} \leq 1 \right\} \\
  & = \| T_n \|.
  \end{split}
\]
Since $m \geq n$, was arbitrary, this means that $\| T \|_\cb = \| T_n \|$.
\end{proof}
\par
Our next lemma is related to \cite[Lemma]{Los} (following \cite{ER},
$\tensor^\lambda$ stands for the injective tensor product of Banach spaces):
\begin{lemma} \label{viktor}
Let $\A$ be a commutative $\cstar$-algebra, and let $n \in \posints$. Then
the canonical map from $M_n \tensor^\lambda M_n(\A)$ to $M_{n^2}(\A)$
has norm at most $n$.
\end{lemma}
\begin{proof}
Again, suppose that $\A \cong {\cal C}_0(\Omega)$ for some locally compact
Hausdorff space $\Omega$. 
\par
We may identify $M_{n^2}(\A)$ with ${\cal C}_0\left(\Omega, M_{n^2}\right)$. 
It is sufficient to show that, for each $\omega \in \Omega$, the map
\begin{equation} \label{yum}
  M_n \tensor^\lambda M_n(\A) \to M_{n^2}, \quad \alpha \tensor f 
  \mapsto \alpha \tensor f(\omega)
\end{equation}
has norm at most $n$.
\par
Let $\omega \in \Omega$, and note that (\ref{yum}) is the composition of the
contraction 
\[
  M_n \tensor^\lambda M_n(\A) \to M_n \tensor^\lambda M_n, 
  \quad \alpha \tensor f 
  \mapsto \alpha \tensor f(\omega)
\]
with the canonical map from $M_n \tensor^\lambda M_n \to M_{n^2}$, which has
norm $n$ by \cite[Lemma]{Los}. Hence, (\ref{yum}) has norm $n$.
\end{proof}
\begin{lemma} \label{norms}
Let $E$ be an operator space, let $\A$ be a commutative $\cstar$-algebra,
and let $n \in \posints$. Then every bounded linear map $T \!: E \to M_n(\A)$
is completely bounded such that $\| T \|_\cb \leq n \| T \|$.
\end{lemma}
\begin{proof}
We can suppose without loss of generality that $E$ is a minimal operator
space, so that, in particular, $M_m(E) = M_m \tensor^\lambda E$ for all $m
\in \posints$.
\par
By Lemma \ref{smithie}, it is enough to show that $\| T_n \| \leq n \| T \|$.
The map $T_n \!: M_n \tensor^\lambda E \to M_{n^2}(\A)$, however, is the 
composition of $\id_{M_n} \tensor T : M_n \tensor^\lambda E \to M_n
\tensor^\lambda M_n(\A)$, which has the same norm as $T$, and the canonical
map from $M_n \tensor^\lambda M_n(\A)$ to $M_{n^2}(\A)$, which has norm at most
$n$ by Lemma \ref{viktor}. Hence, $\| T_n \|$ has norm at most $n \| T \|$.
\end{proof}
\begin{corollary} \label{oscor}
Let $E$ be an operator space, let $\A_1, \ldots, \A_k$ be commutative
$\cstar$-algebras, let $n_1, \ldots, n_k \in \posints$, and let
\[
  \A = M_{n_1}(\A_1) \oplus_\infty \cdots \oplus_\infty M_{n_k}(\A_k).
\] 
Then every bounded linear map $T \!: E \to \A$ is completely bounded such
that $\| T \|_\cb \leq \max \{ n_1, \ldots, n_k \} \| T \|$.
\end{corollary}
\begin{proof}
For $j =1, \ldots, n$, let $T_j \!: E \to M_{n_j}(\A_j)$ be the composition
of $T$ with the projection from $\A$ onto $M_{n_j}(\A_j)$. It follows that
\[
  \begin{split}
  \| T \|_\cb & = \max \{ \| T_1 \|_\cb, \ldots, \| T_k \|_\cb \} \\
  & \leq \max \{ n_1 \| T_1 \|, \ldots, n_k \| T_k \| \},
  \qquad\text{by Lemma \ref{norms}}, \\
  & \leq \max \{ n_1, \ldots, n_k \}
  \max \{\| T_1 \|, \ldots, \| T_k \| \} \\
  & =  \max \{ n_1, \ldots, n_k \} \| T \|,
  \end{split}
\]
which proves the claim.
\end{proof}
\par
As in \cite{ER}, we write $\Tensor$ for projective tensor product of operator
spaces (as opposed to $\tensor^\gamma$). Given two operator spaces $E$ and $F$, 
we have a canonical
contraction from $E \tensor^\gamma F$ to $E \Tensor F$, and, generally, this 
is all that can be said about the relation between $E \tensor^\gamma F$ 
and $E \Tensor F$.
\par
In special situations, however, stronger statements can be made:
\begin{proposition} \label{degG}
Let $G$ be a locally compact group such that $\deg(G) < \infty$, and let $E$
be an operator space. Then the canonical contraction from $A(G) \tensor^\gamma
E$ into $A(G) \Tensor E$ is a topological isomorphism whose inverse has
norm at most $\deg(G)$.
\end{proposition}
\begin{proof}
Let $\VN(G)$ denote the group von Neumann algebra of $G$, and recall that
$A(G)^\ast = \VN(G)$.
\par
We approach the problem from a dual point of view, and show that every
bounded linear map $T \!: E \to \VN(G)$ is completely bounded with
$\| T \|_\cb \leq \deg(G) \| T \|$.
\par
Since $\deg(G) < \infty$, basic structure theory for von Neumann algebras
yields that there are commutative von Neumann algebras $\M_1, \ldots, \M_k$
as well as $n_1, \ldots, n_k \in \posints$ --- with 
$\max \{ n_1, \ldots, n_k \} \leq \deg(G)$ --- such that
\[
  \VN(G) \cong M_{n_1}(\M_1) \oplus_\infty \cdots \oplus_\infty
  M_{n_k}(\M_k).
\]
\par
By Corollary \ref{oscor}, we have a canonical --- obviously $w^\ast$-$w^\ast$-
continuous --- bijection from ${\cal B}(E,\VN(G))$ to $\CB(E,\VN(G))$ with 
norm not exceeding $\max \{ n_1, \ldots, n_k \} \leq \deg(G)$. It follows
that the preadjoint from $A(G) \Tensor E$ to $A(G) \tensor^\gamma E$
of this map, which is the identity on $A(G) \tensor E$, also has norm at 
most $\deg(G)$.
\end{proof}
\begin{remark}
By \cite{Tho} or \cite{Moo}, $\deg(G) < \infty$ holds if and only if $G$ is 
almost abelian. Hence, what we actually show in the proof of Proposition 
\ref{degG}, is that ${\cal B}(A(G),E) = \CB(A(G),E)$ --- not necessarily
with identical norms --- for every almost abelian, locally compact group $G$ 
and every operator space $E$: this result was already proven by Forrest and 
P.\ J.\ Wood (\cite[Theorem 4.5]{FW}). Our approach, however, yields better
norm estimates. If $G$ is a locally compact group, $H$ is a (closed) 
abelian subgroup of $G$ with finite index, $E$ is any operator space, 
and $T \!: A(G) \to E$ is a bounded, linear operator, then an inspection of 
the proof of \cite[Theorem 4.5]{FW} shows that $\| T \|_\cb \leq 
[G : H]  \| T \|$. Proposition \ref{degG}, on the other hand, 
yields the estimate $\| T \|_\cb \leq \deg(G) \| T \|$. In view
of Proposition \ref{thoma} below and the example following it, this
latter estimate is better.
\end{remark}
\begin{corollary} \label{vikcor} 
Let $G$ and $H$ be locally compact groups such that $\deg(G) < \infty$. 
Then the canonical contraction from $A(G) \tensor^\gamma A(H)$ into 
$A(G) \Tensor A(H)$ is a topological isomorphism whose inverse has
norm at most $\deg(G)$.
\end{corollary}
\par
We can now prove the main result of this section:
\begin{theorem} \label{above}
Let $G$ be a locally compact group. Then $\AM_{A(G)} \leq \deg(G)$ holds.
\end{theorem}
\begin{proof}
Since the claim is trivial if $\deg(G) = \infty$, we can suppose that
$\deg(G) < \infty$. Then $G$ is, in particular, amenable. By
\cite[Theorem 3.6]{Rua}, this means that $A(G)$ is {\it operator amenable\/},
i.e.\ there is a bounded net $( d_\alpha )_\alpha$ in $A(G) \Tensor A(G)$
such that (\ref{diag1}) and (\ref{diag2}) hold (with $\Tensor$ instead
of $\tensor^\gamma$). An inspection of the proof of \cite[Theorem 3.6]{Rua}
shows that $( d_\alpha )_\alpha$ can be chosen to have bound one. By Corollary
\ref{vikcor}, $( d_\alpha )_\alpha$ can be viewed as a net in 
$A(G) \tensor^\gamma A(G)$, bounded by $\deg(G)$. Hence, $A(G)$ is
$\deg(G)$-amenable.
\end{proof}
\par
Let $G$ a locally compact group, and let $H$ be a closed, abelian subgroup of 
$G$ with finite index. Then $A(G)$ is amenable by \cite[Theorem 4.1]{LLW},
and an inspection of the proof of \cite[Theorem 4.1]{LLW} shows that 
$A(G)$ is, in fact, $[G:H]$-amenable. 
\par
We shall devote the remainder of this section to showing that Theorem 
\ref{above} provides a better estimate.
\par
The following was proved in \cite{Tho} for the case of a {\it normal\/}
subgroup (\cite[Satz 5]{Tho}):
\begin{proposition} \label{thoma}
Let $G$ be a group, and let $H$ be an abelian subgroup of $G$ of finite index.
Then $\deg(G) \leq [G:H]$ holds.
\end{proposition}
\begin{proof}
Set $n := [G:H]$, and let $\pi \in \hat{G}$. It is well known 
(\cite{Tho} or \cite{Moo}) that $m := \deg(\pi) < \infty$.
\par
We may view $\pi$ as a $^\ast$-representation of the Banach $^\ast$-algebra
$\ell^1(G)$. Then $\pi(\ell^1(G))$ is isomorphic to the $\cstar$-algebra
$M_m$ of all complex $m \times m$-matrices and $\pi(\ell^1(H))$ is a
commutative $\cstar$-subalgebra of $M_m$. This commutative $\cstar$-algebra
is contained in a maximal commutative $\cstar$-subalgebra of $M_m$, and
since --- up to unitary equivalence --- there is only one such
$\cstar$-subalgebra of $M_m$, namely the diagonal matrices, we conclude that
$\dim \pi(\ell^1(H)) \leq m$. 
\par
Let $x_1, \ldots, x_n \in G$ be representatives of the left cosets of $H$,
and note that
\[
  \dim \pi(\ell^1(x_j H)) = \dim \pi(x_j)\pi(\ell^1(H)) = 
  \dim \pi(\ell^1(H)) \leq m \qquad (j =1, \ldots, n).
\]
Since $\ell^1(G) = \ell^1(x_1 H) \oplus \cdots \oplus \ell^1(x_n H)$, we
conclude that
\[
  m^2 = \dim M_m = \dim \pi(\ell^1(G))
  \leq \sum_{j=1}^n \dim \pi(\ell^1(x_j H)) \leq nm,
\]
so that $m \leq n$.
\end{proof}
\par
The inequality in Proposition \ref{thoma} can be strict as the following
example shows:
\begin{example}
Let $\mathrm{A}_5$ be the alternating group in five symbols, i.e.\ the group of all 
even  permutations of $\{1, \ldots, 5 \}$. According to \cite{Con}, 
$\hat{\mathrm{A}}_5$ 
consists of five elements, $\pi_1, \ldots, \pi_5$ say, with
\[
  \deg(\pi_1) = 1, \qquad \deg(\pi_2) = \deg(\pi_3) = 3, \qquad
  \deg(\pi_4) = 4, \qquad\text{and}\qquad \deg(\pi_5) =5,
\]
so that
\[
  \AM_{A(\mathrm{A}_5)} = \frac{61}{15} = 4.0666\ldots \leq 5 = \deg(\mathrm{A}_5). 
\]
Let $H$ be an abelian subgroup of $\mathrm{A}_5$. Assume that $[\mathrm{A}_5:H] = 5$.
Then $H$ is contained in a maximal subgroup, $M$ say, of $\mathrm{A}_5$ whose index
is necessarily at most $5$. Again according to \cite{Con}, $\mathrm{A}_5$ has --- up to 
conjugacy --- only three maximal subgroups whose indices are $5$, $6$, and 
$10$, respectively, so that $[\mathrm{A}_5:M]=5$ and thus $M = H$. The (up to
conjugacy) unique subgroup of $\mathrm{A}_5$ with index $5$, however, is
isomorphic to $\mathrm{A}_4$, the alternating group in four symbols, and therefore not abelian.
It follows that $[\mathrm{A}_5:H]>5$.
\end{example}
\section{The case $\AM_{A(G)} = 1$}
Let $G$ be a locally compact, almost abelian group. Then Theorem \ref{above}
provides us with an estimate for $\AM_{A(G)}$ from above. In view of
(\ref{Barry}) it is clear that it would be naive to expect a similarly simple
estimate from below. 
\par
Nevertheless, some sort of estimate from below is possible.
\par
By $B(G)$, we denote the {\it Fourier--Stieltjes algebra\/} of a locally
compact group $G$ (see \cite{Eym}). For any locally compact group $G$, we use
$G_d$ to denote the same group equipped with the discrete topology.
Finally, the {\it anti-diagonal\/} of a group $G$ is the subset
\[
  G_\Gamma := \{ (x,x^{-1}) : x \in G \}
\]
of $G \times G$, whose indicator function we denote
$\chi_\Gamma$. 
\begin{lemma} \label{below}
Let $G$ be a locally compact group such that $A(G)$ is amenable. Then 
the $\chi_\Gamma$ lies in $B(G_d \times G_d)$ and satisfies 
$\| \chi_\Gamma \|_{B(G_d \times G_d)} \leq \AM_{A(G)}$.
\end{lemma}
\begin{proof}
For any function $f \!: G \to \comps$, define 
\[
  \check{f} \!: G \to G, \quad x \mapsto f(x^{-1}).
\]
The map
\[
  \mbox{}^\vee \!: A(G) \to A(G), \quad f \mapsto \check{f}
\]
is an isometry (see \cite{Eym}).
\par
Let $( d_\alpha )_\alpha$ be an approximate diagonal for $A(G)$ bounded
by $\AM_{A(G)}$. Since $^\vee \!: A(G) \to A(G)$ is an isometry, 
$( (\id \tensor \mbox{}^\vee\ ) d_\alpha )_\alpha$
is a net in $A(G) \tensor^\gamma A(G)$ that is also bounded by 
$\AM_{A(G)}$. We have a canonical contraction from $A(G) \tensor^\gamma A(G)$
into $B(G_d \tensor G_d)$ and may thus view $( (\id \tensor \mbox{}^\vee\ )
d_\alpha )_\alpha$ as a net in $B(G_d \times G_d)$ bounded by $\AM_{A(G)}$. 
From (\ref{diag1}) and (\ref{diag2}), it follows that
$( (\id \tensor \mbox{}^\vee\ ) d_\alpha )_\alpha$ converges to
$\chi_\Gamma$ pointwise on $G \times G$. By \cite[(2.25) Corollaire]{Eym}, 
this means that $\chi_\Gamma \in B(G_d \times G_d)$ with 
$\| \chi_\Gamma \|_{B(G_d \times G_d)} \leq \AM_{A(G)}$.
\end{proof}
\par
Lemma \ref{below} can be used to give a more direct proof of
\cite[Theorem 2.3]{FR}. 
\par
Recall that the {\it coset ring\/} $\Omega(G)$ of a group $G$ is the ring of 
subsets of $G$ generated by all left cosets of subgroups of $G$.
\begin{proposition} \label{almostab}
The following are equivalent for a group $G$:
\begin{items}
\item $G$ is almost abelian;
\item $G_\Gamma \in \Omega(G \times G)$;
\item $\chi_\Gamma \in B(G \times G)$.
\end{items}
\end{proposition}
\begin{proof}
(i) $\Longleftrightarrow$ (ii) is \cite[Proposition 2.2]{FR} and
(ii) $\Longleftrightarrow$ (iii) follows from Host's idempotent theorem
(\cite{Hos}).
\end{proof}
\par
Combining Lemma \ref{below} and Proposition \ref{almostab}, we immediately
recover \cite[Theorem 2.3]{FR}:
\begin{corollary} \label{Brian+I}
The following are equivalent for a locally compact group $G$:
\begin{items}
\item $G$ is almost abelian;
\item $A(G)$ is amenable.
\end{items}
\end{corollary}
\begin{remark}
The present proof for Corollary \ref{Brian+I} is more direct than the one
give in \cite{FR} because it invokes Host's idempotent theorem directly
instead of making the detour over \cite{FKLS}.
\end{remark}
\par
It remains to be seen whether or not Lemma \ref{below} will eventually lead
to a more satisfactory bound from below for the amenability constant of a
Fourier algebra: very little seems to be known on the norms of idempotents
in Fourier--Stieltjes algebras (see the remark below, following 
Theorem \ref{oneam}).
\par
Let $G$ be an abelian locally compact group, so that $A(G) \cong
L^1(\hat{G})$. In view of \cite[Corollary 1.11]{Sto}, this means that 
$\AM_{A(G)} = 1$.
Concluding this note, we shall now see that the locally compact groups $G$
for which $\AM_{A(G)} =1$ are precisely the abelian ones. The two ingredients
of the proof are Lemma \ref{below} and the following proposition that parallels
Proposition \ref{almostab}:
\begin{proposition} \label{ab}
The following are equivalent for a group $G$:
\begin{items}
\item $G$ is abelian;
\item $G_\Gamma$ is a subgroup of $G \times G$;
\item $\chi_\Gamma$ lies in $B(G \times G)$ and has norm one.
\end{items}
\end{proposition}
\begin{proof}
(i) $\Longleftrightarrow$ (ii) is straightforward.
\par
(ii) $\Longrightarrow$ (iii): If $G_\Gamma$ is a subgroup of $G \times G$,
its indicator function is positive definite so that 
$\| \chi_\Gamma \|_{B(G \times G)} = \chi_\Gamma(e,e) = 1$.
\par
(iii) $\Longrightarrow$ (ii): By \cite[Theorem 2.1]{IS},
$G_\Gamma$ must be a left coset of some subgroup of $G \times G$. Since
$(e,e) \in G_\Gamma$, this means that $G_\Gamma$ is, in fact, a subgroup of
$G \times G$.
\end{proof}
\begin{theorem} \label{oneam}
The following are equivalent for a locally compact group $G$:
\begin{items}         
\item $G$ is abelian;
\item $\AM_{A(G)} = 1$.
\end{items}
\end{theorem}
\begin{proof}
We have already observed that (i) $\Longrightarrow$ (ii) holds. The converse
is an immediate consequence of Lemma \ref{below} and Proposition \ref{ab}.
\end{proof}
\begin{remark}
If $G$ is a finite, non-abelian group, $\AM_{A(G)} \geq \frac{3}{2}$ holds
by \cite[Proposition 4.3]{Joh3}. 
In view of Theorem \ref{oneam}, one wonders if this 
estimate from below still holds for arbitrary locally compact groups (with
$\frac{3}{2}$ possibly replaced by another universal constant strictly greater
than one). In view
of Lemma \ref{below}, one way of obtaining such a constant would be to find
an estimate for $\| \chi_\Gamma \|$ from below. In \cite{Sae}, it is proved
for {\it abelian\/} $G$ that, the norm of an idempotent in $B(G)$ is either 
one or at least $\frac{1}{2}\left( 1+ \sqrt{2}\right)$. A similar dichotomy 
result for general, locally compact groups $G$ would immediately yield
a universal bound (strictly greater than one) from below for $\AM_{A(G)}$
for non-abelian $G$.
\end{remark}
\renewcommand{\baselinestretch}{1.0}
\dated
\vfill
\renewcommand{\baselinestretch}{1.2}
\begin{tabbing}
{\it Author's address\/}: \= Department of Mathematical and Statistical Sciences \\
                 \> University of Alberta \\
                 \> Edmonton, Alberta \\
                 \> Canada, T6G 2G1 \\[\medskipamount]
{\it E-mail\/}:  \> {\tt vrunde@ualberta.ca} \\[\medskipamount]
{\it URL\/}: \> {\tt http://www.math.ualberta.ca/$^\sim$runde/}
\end{tabbing}                                                                           
\end{document}